\documentclass[research]{elsarticle}

\usepackage{lineno,hyperref}
\usepackage{amssymb}
\usepackage{graphicx}
\usepackage{amsmath}
\usepackage{amsfonts,color}
\usepackage{amsfonts}%
\usepackage{float}
\usepackage{multirow}
\usepackage{diagbox}

\modulolinenumbers[5]

\journal{Applied Mathematical Modelling}

\providecommand{\U}[1]{\protect\rule{.1in}{.1in}}
\providecommand{\U}[1]{\protect\rule{.1in}{.1in}}
\providecommand{\U}[1]{\protect\rule{.1in}{.1in}}
\def\barr{\begin{array}}
\def\earr{\end{array}}
\def\bec#1{\begin{equation}\label{#1}}
\def\becn{\begin{equation*}}
\def\endec{\end{equation}}
\def\endecn{\end{equation*}}

\def\vi{\raise 2pt\hbox{,}}

\newtheorem{theorem}{Theorem}

\newtheorem{remark}[theorem]{Remark}

\DeclareGraphicsExtensions{.eps,.bmp,.jpg,.pdf,.mps,.png,.gif}

\begin{document}

\begin{frontmatter}

\title{\textbf{Qualitative properties of solutions in the time differential dual-phase-lag model of heat conduction}}

\author[mymainaddress]{Stan Chiri\c t\u a \corref{mycorrespondingauthor} }
\cortext[mycorrespondingauthor]{Corresponding author}
\ead{schirita@uaic.ro}

\author[mysecondaryaddress]{Michele Ciarletta }
\ead{mciarletta@unisa.it}

\author[mysecondaryaddress]{Vincenzo Tibullo }
\ead{vtibullo@unisa.it}

\address[mymainaddress]{Faculty of Mathematics, Al. I. Cuza University of Ia\c{s}i,
700506 -- Ia\c{s}i,\\ \&  Octav Mayer Mathematics Institute, Romanian Academy,
700505 -- Ia\c si, Romania}

\address[mysecondaryaddress]{Dipartimento di Ingegneria Industriale/DIIN, University of Salerno,  84084 - Fisciano (SA), Italia}

\address[mysecondaryaddress]{Dipartimento di Matematica, University of Salerno, 84084 - Fisciano (SA), Italia}

\setcounter{theorem}{0}

\begin{abstract}
\noindent In this paper we study the time differential dual-phase-lag model of heat conduction incorporating the microstructural interaction effect in the fast-transient process of heat transport. We analyse the influence of the delay times upon some qualitative properties of the solutions of the initial boundary value problems associated to such a model. Thus, the uniqueness results are established under the assumption that the conductivity tensor is positive definite and the delay times $\tau_q$ and $\tau_T$ vary in the set $\{0\leq \tau_q\leq 2\tau_T\}\cup \{0<2\tau_T< \tau_q\}$. For the continuous dependence problem we establish two different estimates. The first one is obtained for the delay times with $0\leq \tau_q \leq 2\tau_T$, which agrees with the thermodynamic restrictions on the model in concern, and the solutions are stable. The second estimate is established for the delay times with $0<2\tau_T< \tau_q$ and it allows the solutions to have an exponential growth in time. The spatial behavior of the transient solutions and the steady-state vibrations is also addressed. For the transient solutions we establish a theorem of influence domain, under the assumption that the delay times are in $\left\{0<\tau_q\leq 2\tau_T\right\}\cup \left\{0<2\tau_T<\tau_q\right\}$. While for the amplitude of the harmonic vibrations we obtain an exponential decay estimate of Saint-Venant type, provided the frequency of vibration is lower than a critical value and without any restrictions upon the delay times.

Link ot publisher: \url{https://doi.org/10.1016/j.apm.2017.05.023}
\bigskip

\noindent\textit{Key Words:} dual-phase-lag heat conduction model, uniqueness, continuous dependence, spatial behavior, transient solutions, steady-state vibrations

\end{abstract}

\end{frontmatter}

\section*{Introduction}
\label{intro}
In this paper we study the dual-phase-lag model of heat conduction incorporating the microstructural interaction effects in the fast-transient process of heat transport. Tzou \cite{Tzou1995a,Tzou1995b,Tzou1995c} (see also \cite{Tzou2015} and \cite{Chandra1998} and the references therein) proposed the following time differential constitutive law for the heat flux vector $q_i$
\begin{equation}
\begin{split}
&q_i\left(\mathbf{x},t\right)+\tau_q \frac{\partial q_i}{\partial t}\left(\mathbf{x},t\right)+\frac12 \tau_q^2 \frac{\partial^2 q_i}{\partial t^2}\left(\mathbf{x},t\right)\\[4mm]
&=-k_{ij}(\mathbf{x})T_{,j}\left(\mathbf{x},t\right)-\tau_{T} k_{ij}(\mathbf{x})\frac{\partial  T_{,j}}{\partial t}\left(\mathbf{x},t\right),\label{1.1}
\end{split}
\end{equation}
where $\tau_q \geq 0$ and $\tau_T \geq 0$ are the delay times, $T$ is the temperature variation and $k_{ij}$ is the conductivity tensor. It was established by Fabrizio and Lazzari \cite{Fabrizio2014} that the restrictions imposed by thermodynamics on the constitutive equation (\ref{1.1}), within the framework of a linear rigid conductor, implies that the delay times have to satisfy the inequality $0\leq \tau_q \leq 2\tau_T$.

When the constitutive equation (\ref{1.1}) is coupled with the energy equation
\begin{equation}
-q_{i,i}+\varrho r=a\frac{\partial T}{\partial t},\label{1.2}
\end{equation}
and under appropriate regularity assumptions, then we obtain the following governing equation of hyperbolic type for temperature field $T$
\begin{equation}
\begin{split}
&\left(1+\tau_q\, \frac{\partial}{\partial t}+\frac12 \, \tau_q^2\, \frac{\partial^2}{\partial t^2}\right)\left(a \frac{\partial T}{\partial t} -\varrho r\right)\\[4mm]
&=\left[k_{ij}\left(1+\tau_T\, \frac{\partial}{\partial t}\right)T_{,j}\right]_{,i}.\label{1.3}
\end{split}
\end{equation}
Such equation was studied intensively in literature in many papers (see, for example,
\cite{Wang2001}-\cite{Abouelregal2015}).
 Quintanilla \cite{Quinta2002} has shown that the equation (\ref{1.3}), together with suitable initial conditions for $T$, $\left(\partial T\right)/\left(\partial t\right)$ and $\left(\partial^2 T\right)/\left(\partial t^2\right)$ and appropriate boundary conditions in terms of $T$, leads to an exponentially stable system when $0< \tau_q < 2\tau_T$ and to an unstable system when $0< 2\tau_T < \tau_q$. While in
\cite{Wang2001}-\cite{Xu2002}
the well-posedness problem is studied, provided some appropriate restrictions upon the parameters $\tau_q$ and $\tau_T$ are assumed. At this regard, Quintanilla \cite{Quinta2008,Quinta2009} performed a very interesting analysis about the well-posed problems of dual- and three-phase-lag models of heat conduction equation. We have to outline that some approximations of the Tzou's theory have been studied recently by Quintanilla \cite{Quinta2016}  and Amendola \emph{et al.} \cite{Amendola2016}.

In \cite{Chirita2015} the authors investigate the propagation of plane time harmonic waves and surface waves in the case of a homogeneous and isotropic dual-phase-lag thermoelastic material.

In the present paper we formulate the initial boundary value problem of heat conduction model based on the constitutive equation (\ref{1.1}) and the basic energy equation (\ref{1.2}), considering it is described by a differential system for the unknown couple $\{T,q_i\} $. That means we consider initial conditions for $T$, $q_i$ and $\left(\partial q_i\right)/\left(\partial t\right)$ and appropriate boundary conditions in terms of temperature variation and heat flux vector. Then we study the effects of the presence of the delay times upon the uniqueness and continuous data dependence results as well as the spatial behavior of the transient solutions and of the harmonic vibrations. The uniqueness results are established without any restrictions upon the delay times, except the case of the class of materials characterized by zero delay time of phase lag of the conductive temperature gradient and for which the delay time in the phase lag of heat flux vector is strictly positive. In such a case it should be expected to have an ill-posed model. Moreover, it was shown by Fabrizio and Lazzari \cite{Fabrizio2014} that the corresponding model (with $\tau_T=0$ and $\tau_q>0$) is incompatible with the thermodynamic principles.

We also address the problem of continuous dependence of solutions with respect to the given data. To this aim we establish appropriate conservation laws and then we use the Gronwall's inequality in order to establish two estimates describing the continuous dependence of solutions with respect to the prescribed initial data and with respect to the given supply terms. The first estimate is established for delay times satisfying the inequality $0\leq \tau_q \leq 2\tau_T$, which is in accord with the thermodynamic restriction established by Fabrizio and Lazzari \cite{Fabrizio2014}. The second estimate is established under the assumption $0<2\tau_T < \tau_q$ and it allows the solutions to have a growth exponential in time.

Finally we proceed to study the spatial behavior of transient solutions and to this end we establish a domain of influence theorem, provided the delay times are in the set $\{0<\tau_q\leq 2\tau_T\}\cup \{0<2\tau_T<\tau_q\}$. Explicit estimates are given, in terms of the constitutive thermal coefficients and delay times, for the speed of signal propagation for each of the two subsets of delay times. Moreover, for the steady-state vibrations we establish an exponential decaying estimate in terms of the amplitude vibration, provided the frequency is lower than a critical value and for delay times $\tau_q\geq 0$ and $\tau_T\geq 0$.

\begin{table}
\begin{tabular}{| c | m{4.5cm} | m{4.7cm} | c}
	\hline
	\diagbox{\parbox[c]{2cm}{Type of \\ solutions \\[-0.3cm]}}{Delay times} 
		& \centering $0\leq\frac{1}{2}\tau_q\leq\tau_T$ 
		& \parbox[c]{4.7cm}{\centering $0<\tau_T<\frac{1}{2}\tau_q$} \\ 
	\hline
	\multirow{3}{*}{\parbox[c][4.0cm][c]{2cm}{Transient \\ 
									solutions}} 
		& \parbox[c][2cm][c]{4.5cm}{\centering Uniqueness results}
		& \parbox[c][2cm][c]{4.7cm}{\centering Uniqueness results} \\[1cm]
	\cline{2-3}
		& \parbox[c][2cm][c]{4.5cm}{Continuous dependence \\
									estimates for models \\
									compatible with \\
									thermodynamics} & 
		  \parbox[c][2cm][c]{4.7cm}{Continuous dependence \\
		  							estimates for solutions \\
									with an exponential \\
									growth in time} \\ 
	\cline{2-3}
    	& \parbox[c][2cm][c]{4.5cm}{Theorem of influence \\ 
									domain, provided $\tau_q>0$} 
		& \parbox[c][2cm][c]{4.7cm}{Theorem of influence \\
									domain } \\ 
	\hline
	\parbox[c][2cm][c]{2cm}{Steady \\ 
					state \\
					solutions} 
		& \multicolumn{2}{|c|}{\parbox[c]{9cm}{Exponential decay estimate of Saint-Venant type, \\
												provided the frequency of vibration is lower \\
												than a critical value}} \\[1cm]
	\hline
\end{tabular}
\caption{Summary of results}
\end{table}

\section*{Formulation of the initial boundary value problem}

We assume that a regular region $B$ is filled by an inhomogeneous and anisotropic material with dual phase lag times.

Throughout this paper we consider the initial boundary value problem $\mathcal{P}$ defined by the field equations (\ref{1.1}) and (\ref{1.2}), the initial conditions
\begin{equation}
\begin{split}
&T(\mathbf{x},0)=T^{0}(\mathbf{x}),\\[4mm]
&q_i(\mathbf{x},0)=q_i^0(\mathbf{x}),\quad \frac{\partial q_i}{\partial t}(\mathbf{x},0)=\dot q_i^0(\mathbf{x}),\quad \textit{on} \quad \overline{B},\label{2.1}
\end{split}
\end{equation}
and the boundary conditions
\begin{equation}
\begin{split}
&T(\mathbf{x},t)=\vartheta(\mathbf{x},t)\quad \textit{on}\quad \overline{\Sigma}_1\times [0,\infty),\\[4mm]
&q_i(\mathbf{x},t)n_i=\xi(\mathbf{x},t)\quad \textit{on}\quad \Sigma_2\times [0,\infty).\label{2.2}
\end{split}
\end{equation}
Here $T^{0}\left(  \mathbf{x}\right)  $, $q_{i}^{0}\left(
\mathbf{x}\right)  $, $\dot q_i^0\left(  \mathbf{x}\right)$ and $ \vartheta \left(  \mathbf{x}%
,t\right)  $, $ \xi \left(  \mathbf{x},t\right)  $ are
prescribed smooth functions. Moreover, $\Sigma_{1}$ and $\Sigma_{2}$ are
subsets of the boundary $\partial B$ so that $\overline{\Sigma}_{1}%
\cup\Sigma_{2}=\partial B$ and $\Sigma_{1}\cap\Sigma_{2}=\emptyset$.

By a solution of the initial boundary value problem $\mathcal{P}$ corresponding to the given data $\mathcal{D}=\{r;T^0,q_i^0,\dot q_i^0;\vartheta,\xi\}$ we mean the ordered array $\mathcal{S}=\{T,q_i\}$ defined on $\overline{B}\times[0,\infty)$ with the properties that $T(\mathbf{x},t)\in C^{1,1}$ $\left(B\times(0,\infty)\right)$, $q_{i}(\mathbf{x},t)\in C^{1,2}\left(B\times(0,\infty)\right)$ and which satisfy the field equations (\ref{1.1}) and (\ref{1.2}), the initial conditions (\ref{2.1}) and the boundary conditions (\ref{2.2}). In what follows we denote by $\mathcal{P}_0$ the initial boundary value problem $\mathcal{P}$ corresponding to the zero given data $\mathcal{D}=\{r;T^0,q_i^0,\dot q_i^0;\vartheta,\xi\}=0$.

For further convenience we introduce the following operators
\begin{equation}
f'(t)=\int_0^tf(s)ds,\quad f''(t)=\int_0^t\int_0^sf(z)dzds, \quad\textit{etc.};\label{2.3}
\end{equation}
\begin{equation}
\hat{g}(t)=g''(t)+\tau_q g'(t)+\frac12\, \tau_q^2 g(t);\label{2.4}
\end{equation}
\begin{equation}
\tilde {h}(t)=h'(t)+\tau_T h(t).\label{2.5}
\end{equation}

It is worth to note that if $\mathcal{S}=\{T,q_i\}$ is a solution of the initial boundary value problem $\mathcal{P}$ corresponding to the given data $\mathcal{D}=\{r;T^0,q_i^0,\dot q_i^0;\vartheta,\xi\}$, then $\hat{\mathcal{S}}=\{\hat{T},\hat{q}_i\}$ is a solution of the initial boundary value problem $\hat{\mathcal{P}}$ defined by the basic equations
\begin{equation}
a\frac{\partial \hat{T}}{\partial t}=-\hat{q}_{i,i}+R,\label{2.6}
\end{equation}
\begin{equation}
\hat{q}_i(t)=-k_{ij}\left(T''_{,j}(t)+\tau_T T'_{,j}(t)\right)+\vartheta_i(t),\label{2.7}
\end{equation}
the initial conditions
\begin{equation}
\begin{split}
&\hat{T}(\mathbf{x},0)=\frac12\, \tau_q^2\, T^0(\mathbf{x}),\\[3mm]
&\hat{q}_i(\mathbf{x},0)=\frac12\, \tau_q^2q_i^0(\mathbf{x}),\quad \frac{\partial \hat{q}_i}{\partial t}(\mathbf{x},0)=\tau_q q_i^0(\mathbf{x})+\frac12\, \tau_q^2 \dot{q}_i^0(\mathbf{x}),\label{2.8}
\end{split}
\end{equation}
and the boundary conditions
\begin{equation}
\begin{split}
&\hat{T}(\mathbf{x},t)=\hat{\vartheta}(\mathbf{x},t)\quad \textit{on}\quad \overline{\Sigma}_1\times [0,\infty),\\[3mm]
&\hat{q}_i(\mathbf{x},t)n_i=\hat{\xi}(\mathbf{x},t)\quad \textit{on}\quad \Sigma_2\times [0,\infty),\label{2.9}
\end{split}
\end{equation}
where
\begin{equation}
\begin{split}
&R(\mathbf{x},t)=\varrho \hat{r}(\mathbf{x},t)+a(t+\tau_q)T^0(\mathbf{x}),\\[3mm]
&\vartheta_i(\mathbf{x},t)=\left[\tau_T\, k_{ij} T^0_{,j}(\mathbf{x})+\tau_q q_i^0(\mathbf{x})+\frac12\, \tau_q^2\, \dot{q}_i^0(\mathbf{x})\right]t +\frac12\, \tau_q^2\, q_i^0(\mathbf{x}).\label{2.10}
\end{split}
\end{equation}

Suppose now that the conductivity tensor $k_{ij}$ is non-singular and let us denote by $K_{ij}$ its inverse tensor so that
\begin{equation}
k_{ij}K_{jk}=K_{ij}k_{jk}=\delta_{ik}.\label{2.11}
\end{equation}
It is a straightforward task to verify that if $\mathcal{S}=\{T,q_i\}$ is a solution of the initial boundary value problem $\mathcal{P}$ corresponding to the given data $\mathcal{D}=\{r;T^0,q_i^0,\dot q_i^0;\vartheta,\xi\}$, then $\tilde{\mathcal{S}}=\{\tilde{T},\tilde{q}_i\}$ is a solution of the initial boundary value problem $\tilde{\mathcal{P}}$ defined by the basic equations
\begin{equation}
a\frac{\partial \tilde{T}}{\partial t}=-\tilde{q}_{i,i}+R^*,\label{2.12}
\end{equation}
\begin{equation}
\tilde{T}_{,i}(t)=-K_{ij}\left[q'_{j}(t)+\tau_q q_{j}(t)+\frac12\, \tau_q^2\, \dot{q}_j(t)\right]+\varphi_i,\label{2.13}
\end{equation}
the initial conditions
\begin{equation}
\begin{split}
&\tilde{T}(\mathbf{x},0)= \tau_T\, T^0(\mathbf{x}),\\[3mm]
&\tilde{q}_i(\mathbf{x},0)= \tau_T \, q_i^0(\mathbf{x}),\quad \frac{\partial \tilde{q}_i}{\partial t}(\mathbf{x},0)= q_i^0(\mathbf{x})+ \tau_T\, \dot{q}_i^0(\mathbf{x}),\label{2.14}
\end{split}
\end{equation}
and the boundary conditions
\begin{equation}
\begin{split}
&\tilde{T}(\mathbf{x},t)=\tilde{\vartheta}(\mathbf{x},t)\quad \textit{on}\quad \overline{\Sigma}_1\times [0,\infty),\\[3mm]
&\tilde{q}_i(\mathbf{x},t)n_i=\tilde{\xi}(\mathbf{x},t)\quad \textit{on}\quad \Sigma_2\times [0,\infty),\label{2.15}
\end{split}
\end{equation}
where
\begin{equation}
\begin{split}
&R^*(\mathbf{x},t)=\varrho \tilde{r}(\mathbf{x},t)+aT^0(\mathbf{x}),\\[3mm]
&\varphi_i(\mathbf{x})=\tau_T\, T^0_{,i}(\mathbf{x})+K_{ij}\left[\tau_q q_j^0(\mathbf{x})+\frac12\, \tau_q^2\, \dot{q}_j^0(\mathbf{x})\right].\label{2.16}
\end{split}
\end{equation}

\section*{Uniqueness results}

Throughout this section we consider the following constitutive hypotheses
\begin{description}
\begin{item}
\item (H1) $k_{ij}$ is a positive definite tensor and hence
\begin{equation}
k_m\xi_i\xi_i\leq k_{ij}\xi_i\xi_j\leq k_M\xi_i\xi_i, \label{3.1}
\end{equation}
where $k_m$ and $k_M$ are the lowest and the greatest eigenvalues of $k_{ij}$;
\end{item}
\begin{item}
\item (H2) $a\neq 0$;
\end{item}
\begin{item}
\item (H3) the delay times are such that
\begin{equation}
\{0\leq\frac12 \, \tau_q\leq \tau_T\}\cup \{0<\tau_T<\frac12 \, \tau_q\};\label{3.2}
\end{equation}
\end{item}
\begin{item}
\item (H4) $meas\, \Sigma_1 \neq 0$.
\end{item}
\end{description}


\bigskip

\begin{theorem}\label{Theorem1}
Suppose that the hypotheses $(H1)$ and $(H3)$ hold true and, moreover, we assume that at least one of the hypotheses (H2) or (H4) is fulfilled. Then the initial boundary value problem $\mathcal{P}$ has at most one solution.
\end{theorem}

\textbf{Proof.} In order to prove the uniqueness result it is sufficient to prove that the zero external given data,
that is $\mathcal{D}=\{r;T^0,q_i^0,\dot q_i^0;\vartheta,\xi\}=0$,
implies that the corresponding solution $\mathcal{S}=\{T,q_i\}$ is vanishing on $\overline{B}\times [0,\infty)$. That means we have to prove that the initial boundary value problem $\mathcal{P}_0$  has only the trivial solution.

Thus, we consider here that $\mathcal{S}=\{T,q_i\}$ is a solution of the initial boundary value problem $\mathcal{P}_0$. It follows then that $\hat{\mathcal{S}}=\{\hat{T},\hat{q}_i\}$ is a solution of the initial boundary value problem $\hat{\mathcal{P}}$ with zero given data, denoted in what follows by $\hat{\mathcal{P}}_0$. Then, by means of an integration with respect to time variable of the equation (\ref{2.6}), followed by the use of the zero initial data, it follows that
\begin{equation}
a\hat{T}(t)=-\int_0^t\hat{q}_{i,i}(z)dz.\label{3.3}
\end{equation}

Furthermore, we start with
\begin{equation}
a\hat{T}(t+s)\hat{T}(t-s)-a\hat{T}(t-s)\hat{T}(t+s)=0,\quad \textit{for all}\quad s\in (0,t),\, t\geq 0,\label{3.4}
\end{equation}
which, by means of (\ref{3.3}), implies
\begin{equation}
\int_B\left[-\hat{T}(t+s)\int_0^{t-s}\hat{q}_{i,i}(z)dz+\hat{T}(t-s)\int_0^{t+s}\hat{q}_{i,i}(z)dz\right]dv=0.\label{3.5}
\end{equation}
If we use the divergence theorem and the zero boundary conditions (\ref{2.9}), from (\ref{3.5}) we get
\begin{equation}
\int_B\left[\hat{T}_{,i}(t+s)\int_0^{t-s}\hat{q}_{i}(z)dz-\hat{T}_{,i}(t-s)\int_0^{t+s}\hat{q}_{i}(z)dz\right]dv=0\label{3.6}
\end{equation}
and hence, by replacing $\hat{q}_i$ from the constitutive equation (\ref{2.7}) and $\hat{T}_{,i}$ by (\ref{2.4}), we obtain
\begin{equation}
\begin{split}
&\int_B\left[k_{ij}T'''_{,j}(t-s)T''_{,i}(t+s)-k_{ij}T'''_{,j}(t+s)T''_{,i}(t-s)\right]dv\\[4mm]
&+\tau_q\int_B\left[k_{ij}T'''_{,j}(t-s)T'_{,i}(t+s)-k_{ij}T'''_{,j}(t+s)T'_{,i}(t-s)\right]dv\\[4mm]
&+\frac12 \, \tau_q^2\int_B\left[k_{ij}T'''_{,j}(t-s)T_{,i}(t+s)-k_{ij}T'''_{,j}(t+s)T_{,i}(t-s)\right]dv\\[4mm]
&+\tau_T\tau_q\int_B\left[k_{ij}T''_{,j}(t-s)T'_{,i}(t+s)-k_{ij}T''_{,j}(t+s)T'_{,i}(t-s)\right]dv\\[4mm]
&+\frac12\, \tau_T\tau_q^2\int_B\left[k_{ij}T''_{,j}(t-s)T_{,i}(t+s)-k_{ij}T''_{,j}(t+s)T_{,i}(t-s)\right]dv=0.\label{3.7}
\end{split}
\end{equation}
Moreover, we can write
\begin{equation}
\begin{split}
&\frac{\partial}{\partial s}\left\{\int_B k_{ij}T'''_{,j}(t-s)T'''_{,i}(t+s)dv+\tau_T\tau_q\, \int_B k_{ij}T''_{,j}(t-s)T''_{,i}(t+s)dv\right.\\[4mm]
&+\tau_q\, \int_B\left[k_{ij}T'''_{,j}(t-s)T''_{,i}(t+s)+k_{ij}T'''_{,j}(t+s)T''_{,i}(t-s)\right]dv\\[4mm]
&+\frac12\, \tau_q^2\int_B\left[k_{ij}T'''_{,j}(t-s)T'_{,i}(t+s)+k_{ij}T'''_{,j}(t+s)T'_{,i}(t-s)\right]dv\\[4mm]
&+\frac12\, \tau_T\tau_q^2\int_B\left[k_{ij}T''_{,j}(t-s)T'_{,i}(t+s)+k_{ij}T''_{,j}(t+s)T'_{,i}(t-s)\right]dv\\[4mm]
&\left.+\frac12\, \tau_q^2\int_B k_{ij}T''_{,j}(t-s)T''_{,i}(t+s)dv\right\}=0,\label{3.8}
\end{split}
\end{equation}
so that an integration with respect to $s$ on the interval $[0,t]$ implies
\begin{equation}
\begin{split}
&\frac{d}{dt}\left\{\tau_q\int_B k_{ij}T'''_{,j}(t)T'''_{,i}(t)dv+\tau_q^2\int_B k_{ij}T'''_{,j}(t)T''_{,i}(t)dv\right.\\[4mm]
&\left.+\frac12 \, \tau_T\tau_q^2\int_B k_{ij}T''_{,j}(t)T''_{,i}(t)dv\right\}+\int_B k_{ij}T'''_{,j}(t)T'''_{,i}(t)dv\\[4mm]
&+\tau_q\left(\tau_T-\frac{\tau_q}{2}\right)\int_B k_{ij}T''_{,j}(t)T''_{,i}(t)dv=0.\label{3.9}
\end{split}
\end{equation}

Thus, by twice integrations with respect to time variable, from (\ref{3.9}) we obtain the following identity
\begin{equation}
\begin{split}
&\int_0^t\int_0^s\int_B k_{ij}T'''_{,j}(z)T'''_{,i}(z)dvdzds+\tau_q\int_0^t\int_B k_{ij}T'''_{,j}(s)T'''_{,i}(s)dvds\\[4mm]
&+\frac{\tau_q^2}{2}\int_B k_{ij} T'''_{,j}(t)T'''_{,i}(t)dv+\frac12\, \tau_T\tau_q^2\int_0^t\int_B k_{ij}T''_{,j}(s)T''_{,i}(s)dvds\\[4mm]
&+\tau_q\left(\tau_T-\frac{\tau_q}{2}\right)\int_0^t\int_0^s\int_B k_{ij}T''_{,j}(z)T''_{,i}(z)dvdzds=0,\quad t\geq 0. \label{3.10}
\end{split}
\end{equation}

Let us first suppose that
\begin{equation}
0\leq\tau_q\leq 2\tau_T.\label{3.11}
\end{equation}
Then all the terms of the identity (\ref{3.10}) are positive and hence we deduce
\begin{equation}
k_{ij}T'''_{,j}(t)T'''_{,i}(t)=0,\label{3.12}
\end{equation}
so that we have
\begin{equation}
T'''_{,i}(\mathbf{x},t)=0\quad \textit{for all}\quad (\mathbf{x},t)\in B\times (0,\infty).\label{3.13}
\end{equation}
Thus, we have
\begin{equation}
T_{,i}(\mathbf{x},t)=0\quad \textit{for all}\quad (\mathbf{x},t)\in B\times (0,\infty),\label{3.14}
\end{equation}
which, by using the constitutive equation (\ref{2.7}), gives
\begin{equation}
\hat{q}_i(\mathbf{x},t)=0\quad \textit{for all}\quad (\mathbf{x},t)\in B\times (0,\infty).\label{3.15}
\end{equation}
This last relation together with (\ref{2.4}) gives
\begin{equation}
q_i\left(\mathbf{x},t\right)+\tau_q \frac{\partial q_i}{\partial t}\left(\mathbf{x},t\right)+\frac12 \tau_q^2 \frac{\partial^2 q_i}{\partial t^2}\left(\mathbf{x},t\right)=0,\label{3.16}
\end{equation}
which, with the aid of the initial conditions
\begin{equation}
q_i(\mathbf{x},0)=0,\quad \dot q_i(\mathbf{x},0)=0,\label{3.17}
\end{equation}
implies
\begin{equation}
q_i(\mathbf{x},t)=0\quad \textit{for all}\quad (\mathbf{x},t)\in \overline{B}\times [0,\infty).\label{3.18}
\end{equation}
Furthermore, from the relations (\ref{1.2}) and (\ref{3.18}) we have
\begin{equation}
a\frac{\partial T}{\partial t}=0.\label{3.19}
\end{equation}

Now, in view of the hypothesis (H2), the relation (\ref{3.19}) furnishes
\begin{equation}
\frac{\partial T}{\partial t}=0\label{3.20}
\end{equation}
and hence, by using the zero initial condition $T(\mathbf{x},0)=0$, we obtain
\begin{equation}
T(\mathbf{x},t)=0\quad \textit{for all}\quad (\mathbf{x},t)\in \overline{B}\times [0,\infty),\label{3.21}
\end{equation}
and therefore, the corresponding solution $\mathcal{S}=\{T,q_i\}$ is vanishing on $\overline{B}\times [0,\infty)$.

Otherwise, if hypothesis (H4) holds true, then from relation (\ref{3.14}) we deduce again the relation (\ref{3.21}) and so we have that the corresponding solution $\mathcal{S}=\{T,q_i\}$ is vanishing on $\overline{B}\times [0,\infty)$.

Let us now consider the case
\begin{equation}
0<\tau_T<\frac12 \, \tau_q.\label{3.22}
\end{equation}
Then the identity (\ref{3.10}) implies
\begin{equation}
\frac12\, \tau_T\tau_q^2\int_0^t\int_B k_{ij}T''_{,j}(s)T''_{,i}(s)dvds\leq \tau_q\left(\frac{\tau_q}{2}-\tau_T\right)\int_0^t\int_0^s\int_B k_{ij}T''_{,j}(z)T''_{,i}(z)dvdzds,\label{3.23}
\end{equation}
which leads to the following Gronwall type inequality
\begin{equation}
\Psi(t)\leq \left(\frac{1}{\tau_T}-\frac{2}{\tau_q}\right)\int_0^t\Psi(s)ds,\quad t\geq 0,\label{3.24}
\end{equation}
with
\begin{equation}
\Psi(t)=\int_0^t\int_B k_{ij}T''_{,j}(s)T''_{,i}(s)dvds.\label{3.25}
\end{equation}
Then the Gronwall's lemma furnishes
\begin{equation}
\Psi(t)=0\quad \textit{for all}\quad t\geq 0\label{3.26}
\end{equation}
and hence we have
\begin{equation}
T''_{,i}(\mathbf{x},t)=0\quad \textit{for all}\quad (\mathbf{x},t)\in B\times (0,\infty).\label{3.27}
\end{equation}
This last relation leads to the conclusion expressed in (\ref{3.14}) and the analysis follows the way described in the proof of the case (\ref{3.11}). Thus, the proof is complete.

A similar result can be obtained if we replace the hypothesis $a\neq 0$ by $a>0$. Thus, we have

\begin{theorem}\label{Theorem2}
Suppose that the hypotheses $(H1)$ holds true, $a>0$ and $\tau_T> 0$ and $\tau_q\geq 0$. Then the initial boundary value problem $\mathcal{P}$ has at most one solution.
\end{theorem}

\textbf{Proof.} In order to prove the result we use the time weighted method \cite{ChirCiar1999}. Thus, we use the time weight function $e^{-\alpha t}$, $\alpha\geq 0$ and the basic equations of the initial boundary value problem $\mathcal{P}_0$ to get the following time weighted conservation law
\begin{equation}
\begin{split}
&\frac12\, \int_0^t\int_B\, e^{-\alpha s}\left[a\hat{T}^2(s)+\left(\tau_T+\tau_q+\alpha\tau_q^2\right)k_{ij}T''_{,j}(s)T''_{,i}(s)\right.\\[4mm]
&\left.+\frac12\, \tau_T\tau_q^2k_{ij}T'_{,j}(s)T'_{,i}(s)\right]dvds+\frac14\, \tau_q^2\int_B e^{-\alpha t}k_{ij}T''_{,j}(t)T'_{,i}(t)dv\\[4mm]
&+\int_0^t\int_0^s\int_B\, e^{-\alpha s}\left\{\frac{\alpha}{2}\, a\hat{T}^2(z)+\left[1+\frac{\alpha}{2}\left(\tau_T+\tau_q\right)+\frac{\alpha^2}{4}\, \tau_q^2\right]k_{ij}T''_{,j}(z)T''_{,i}(z)\right.\\[4mm]
&\left.+\left[\frac{\alpha}{4}\, \tau_T\tau_q^2+\tau_q\left(\tau_T-\frac12 \, \tau_q\right)\right]\right\}dvdzds=0, \quad t\geq0.\label{3.28}
\end{split}
\end{equation}
Further, it can be easily see that it is possible to choose the parameter $\alpha\geq 0$ such that
\begin{equation}
\frac{\alpha}{4}\, \tau_T\tau_q^2+\tau_q\left(\tau_T-\frac12 \, \tau_q\right)\geq 0.\label{3.29}
\end{equation}
With this choice we observe that all the integral terms in the above time weighted identity are made positive and we can follow an analysis similar with that used in the proof of Theorem \ref{Theorem1} to show that $\mathcal{S}=\{T,q_i\}$ is vanishing on $\overline{B}\times [0,\infty)$.

\medskip

\begin{remark}
We have established the uniqueness results described in the Theorem \ref{Theorem1} and the Theorem \ref{Theorem2} by means of the problem $\hat{\mathcal{P}}$, but they can be proved as well by means of the problem $\tilde{\mathcal{P}}$. We leave this task to the reader.
\end{remark}

\section*{Continuous dependence results}

In this section we study the effects of the delay times upon the problem of continuous dependence of solutions with respect to the supply term and the initial given data. To this aim we consider $\mathcal{S}=\{T,q_i\}$ be a solution of the initial boundary value problem $\mathcal{P}$ corresponding to the given data $\mathcal{D}=\{r;T^0,q_i^0,\dot q_i^0;0,0\}$ and then it follows that $\tilde{\mathcal{S}}=\{\tilde{T},\tilde{q}_i\}$ is a solution of the initial boundary value problem $\tilde{\mathcal{P}}$ with zero boundary data. At this stage we introduce the following functional
\begin{equation}
\mathcal{E}(t)=\frac12\, \int_0^t\int_B a\tilde{T}^2(s)dvds+\frac12\, \int_0^t\int_0^s\int_B K_{ij}q'_i(z)q'_j(z)dvdzds,\quad t\geq 0\label{4.1}
\end{equation}
and we note that the relations (\ref{2.12})-(\ref{2.15}) furnish the following conservation law
\begin{equation}
\begin{split}
&\mathcal{E}(t)+\frac{\tau_T+\tau_q}{2}\int_0^t\int_B K_{ij}q'_j(s)q'_i(s)dvds\\[4mm]
&+\frac14\, \tau_T\tau_q^2\int_0^t\int_B K_{ij}q_i(s)q_j(s)dvds+\frac14\, \tau_q^2 \int_B K_{ij}q'_j(t)q'_i(t)dv\\[4mm]
&+\frac12\, \int_0^t\int_0^s\int_B K_{ij}q'_i(z)q'_j(z)dvdzds\\[4mm]
&+\tau_q\left(\tau_T-\frac12\, \tau_q\right)\int_0^t\int_0^s\int_B K_{ij}q_i(z)q_j(z)dvdzds\\[4mm]
&=\frac{t}{2}\left[\int_B a\tilde{T}^2(0)dv+\frac12\, \tau_T\tau_q^2\int_B K_{ij}q_i^0q_j^0dv\right]\\[4mm]
&+\int_0^t\int_0^s\int_B R^*(z)\tilde{T}(z)dvdzds+\int_0^t\int_0^s\int_B\varphi_i(\mathbf{x})q'_i(z)dvdzds\\[4mm]
& +\tau_T\int_0^t\int_B\varphi_i(\mathbf{x})q'_i(s)dvds,\quad t\geq 0.\label{4.2}
\end{split}
\end{equation}

On the other side, we have
\begin{equation}
K_m \xi_i\xi_i \leq K_{ij}\xi_i\xi_j\leq K_M \xi_i\xi_i,\quad \textit{for all}\quad \xi_i\in \mathbf{R}, \label{4.3}
\end{equation}
where $K_m(\mathbf{x})$ and $K_M(\mathbf{x})$ are the lowest and the greatest eigenvalues of the positive definite tensor $K_{ij}$.
On the basis of the Cauchy-Schwarz inequality and the arithmetic-geometric mean inequality, we have
\begin{equation}
\begin{split}
\int_0^t\int_0^s\int_B\varphi_i(\mathbf{x})q'_i(z)dvdzds&\leq \frac{t^2}{4}\, \int_B\frac{1}{\varepsilon_1 K_m}\, \varphi_i(\mathbf{x})\varphi_i(\mathbf{x})dv\\
&+\frac{1}{2}\, \int_0^t\int_0^s\int_B \varepsilon_1 K_{ij}q'_i(z)q'_j(z)dvdzds, \quad \textit{for every}\quad \varepsilon_1>0,\label{4.4}
\end{split}
\end{equation}
and
\begin{equation}
\begin{split}
\int_0^t\int_B\varphi_i(\mathbf{x})q'_i(s)dvds&\leq \frac{t}{2}\, \int_B\frac{1}{\varepsilon_2 K_m}\, \varphi_i(\mathbf{x})\varphi_i(\mathbf{x})dv\\
&+\frac{1}{2}\, \int_0^t\int_B \varepsilon_2 K_{ij}q'_i(s)q'_j(s)dvds, \quad \textit{for every}\quad \varepsilon_2>0.\label{4.5}
\end{split}
\end{equation}

Setting $\varepsilon_1=\varepsilon_2=1$ into the inequalities (\ref{4.4}) and (\ref{4.5}) and then using the result into identity (\ref{4.2}), we obtain the following estimate
\begin{equation}
\begin{split}
&\mathcal{E}(t)+\frac14\, \tau_T\tau_q^2\int_0^t\int_B K_{ij}q_i(s)q_j(s)dvds\\[4mm]
&+\tau_q\left(\tau_T-\frac12\, \tau_q\right)\int_0^t\int_0^s\int_B K_{ij}q_i(z)q_j(z)dvdzds\\[4mm]
&\leq \frac{t}{2}\left[\int_B a\tilde{T}^2(0)dv+\frac12\, \tau_T\tau_q^2\int_B K_{ij}q_i^0q_j^0dv+\tau_T\, \int_B\frac{1}{K_{m}}\, \varphi_i\varphi_i dv\right]\\[4mm]
&+\frac{t^2}{4}\int_B\frac{1}{K_m}\, \varphi_i\varphi_i dv+\int_0^t\int_0^s\int_B R^*(z)\tilde{T}(z)dvdzds,\quad t\geq 0.\label{4.6}
\end{split}
\end{equation}

Then, the following theorem describes the continuous dependence of solutions of the problem $\mathcal{P}$ with respect to the supply term and with respect to the given initial data.

\begin{theorem}
Suppose that the constitutive hypotheses (H1) and (H3) hold true and $a>0$. Then the solution $\mathcal{S}=\{T,q_i\}(\mathbf{x},t)$ of the initial boundary value problem $\mathcal{P}$ depends continuously on the given data $\mathcal{D}=\{r;T^0,q_i^0,\dot q_i^0;0,0\}$ for any finite time interval $[0,S]$, $S>0$. More precisely, for every $t\in [0,S]$, we have the following results:
\begin{description}
\item (i) when
\begin{equation}
0\leq \tau_q \leq 2\tau_T,\label{4.7}
\end{equation}
the estimate
\begin{equation}
\begin{split}
&\sqrt{\mathcal{E}(t)}\leq \frac{1}{\sqrt{2}}\, \int_0^tg(s)ds\\
&+\left\{\frac{S}{2}\int_B\left[a\tilde{T}^2(0)+\frac12\, \tau_T\tau_q^2 K_{ij}q_i^0q_j^0+\frac{\tau_T}{K_m}\varphi_i\varphi_i\right]dv+\frac{S^2}{4}\, \int_B\frac{1}{K_m}\varphi_i\varphi_idv\right\}^{1/2},\label{4.8}
\end{split}
\end{equation}
holds true with
\begin{equation}
g(t)=\left[\int_0^t\int_B\frac{1}{a}R^{*2}(s)dvds\right]^{\frac12};\label{4.9}
\end{equation}
while

\item (ii) when
\begin{equation}
0<2\tau_T<\tau_q,\label{4.10}
\end{equation}
we have the estimate
\begin{equation}
\begin{split}
&\sqrt{\mathcal{F}(t)}\leq \frac{1}{\sqrt{2}}\, \int_0^t \exp\left(\sigma^2 (t-s)\right)g(s)ds\\
&+\exp\left(\sigma^2 t\right)\left\{\frac{S}{2}\int_B\left[a\tilde{T}^2(0)+\frac12\, \tau_T\tau_q^2 K_{ij}q_i^0q_j^0+\frac{\tau_T}{K_m}\varphi_i\varphi_i\right]dv+\frac{S^2}{4}\, \int_B\frac{1}{K_m}\varphi_i\varphi_idv\right\}^{1/2},\label{4.11}
\end{split}
\end{equation}
with
\begin{equation}
\mathcal{F}(t)=\mathcal{E}(t)+\frac14\, \tau_T\tau_q^2\, \int_0^t\int_B K_{ij}q_i(s)q_j(s)dvds,\quad \sigma^2=\frac{1}{\tau_T}-\frac{2}{\tau_q}.\label{4.12}
\end{equation}
\end{description}
\end{theorem}

\textbf{Proof.} Since $a>0$ and in view of the constitutive hypothesis (H1) and (\ref{4.3}), it follows that $\mathcal{E}$, as defined by (\ref{4.1}), can be considered as a measure for $\mathcal{S}$ in the sense that $\mathcal{E}(t)\geq 0$ for all $t\geq 0$ and $\mathcal{E}(t)=0$ implies $\mathcal{S}=0$.

Let us consider first the case (i). Then, a consequence of relation (\ref{4.6}) is the following Gronwall type inequality
\begin{equation}
\begin{split}
&\mathcal{E}(t)\leq \int_0^t g(s)\sqrt{2\mathcal{E}(s)}ds+\frac{S^2}{4}\int_B\frac{1}{K_m}\, \varphi_i\varphi_i dv\\
&+\frac{S}{2}\left[\int_B a\tilde{T}^2(0)dv+\frac12\, \tau_T\tau_q^2\int_B K_{ij}q_i^0q_j^0dv+\tau_T\, \int_B\frac{1}{K_m}\, \varphi_i\varphi_i dv\right].\label{4.13}
\end{split}
\end{equation}
To integrate the integral inequality (\ref{4.13}) we set
\begin{equation}
\begin{split}
&\Phi(t)= \left\{\int_0^t g(s)\sqrt{2\mathcal{E}(s)}ds+\frac{S^2}{4}\int_B\frac{1}{K_m}\, \varphi_i\varphi_i dv\right.\\
&\left.+\frac{S}{2}\left[\int_B a\tilde{T}^2(0)dv+\frac12\, \tau_T\tau_q^2\int_B K_{ij}q_i^0q_j^0dv+\tau_T\, \int_B\frac{1}{K_m}\, \varphi_i\varphi_i dv\right]\right\}^{1/2}.\label{4.14}
\end{split}
\end{equation}
and we note that (\ref{4.13}) implies
\begin{equation}
\sqrt{\mathcal{E}(t)}\leq \Phi(t).\label{4.15}
\end{equation}

On the other side, from (\ref{4.14}) we deduce
\begin{equation}
\frac{d\Phi}{dt}(t)\leq \frac{1}{\sqrt{2}}\, g(t),\label{4.16}
\end{equation}
and hence, by an integration, we get
\begin{equation}
\Phi(t)\leq \Phi(0)+\frac{1}{\sqrt{2}}\int_0^tg(s)ds.\label{4.17}
\end{equation}
Concluding, we see that the relations (\ref{4.14}) to (\ref{4.17}) imply the estimate (\ref{4.8}) and the proof of the case (i) is complete.

Let us further consider the case (ii). Then the relations (\ref{4.6}) and (\ref{4.12}) imply
\begin{equation}
\begin{split}
&\mathcal{F}(t)\leq 2\sigma^2\, \int_0^t\left[\frac{\tau_T\tau_q^2}{4}\, \int_0^s\int_B K_{ij}q_i(z)q_j(z)dvdz\right]ds\\
&+ \frac{t}{2}\left[\int_B a\tilde{T}^2(0)dv+\frac12\, \tau_T\tau_q^2\int_B K_{ij}q_i^0q_j^0dv+\tau_T\, \int_B\frac{1}{K_m}\, \varphi_i\varphi_i dv\right]\\[4mm]
&+\frac{t^2}{4}\int_B\frac{1}{K_m}\, \varphi_i\varphi_i dv+\int_0^t\int_0^s\int_B R^*(z)\tilde{T}(z)dvdzds,\quad t\geq 0,\label{4.18}
\end{split}
\end{equation}
and hence we obtain the following Gronwall type inequality
\begin{equation}
\begin{split}
&\mathcal{F}(t)\leq 2\sigma^2\int_0^t\mathcal{F}(s)ds+\int_0^t g(s)\sqrt{2\mathcal{F}(s)}ds+\frac{S^2}{4}\int_B\frac{1}{K_m}\, \varphi_i\varphi_i dv\\
&+\frac{S}{2}\left[\int_B a\tilde{T}^2(0)dv+\frac12\, \tau_T\tau_q^2\int_B K_{ij}q_i^0q_j^0dv+\tau_T\, \int_B\frac{1}{K_m}\, \varphi_i\varphi_i dv\right].\label{4.19}
\end{split}
\end{equation}
Further, if we set
\begin{equation}
\begin{split}
&\Pi(t)= \left\{2\sigma^2\, \int_0^t\mathcal{F}(s)ds+\int_0^t g(s)\sqrt{2\mathcal{F}(s)}ds+\frac{S^2}{4}\int_B\frac{1}{K_m}\, \varphi_i\varphi_i dv\right.\\
&\left.+\frac{S}{2}\left[\int_B a\tilde{T}^2(0)dv+\frac12\, \tau_T\tau_q^2\int_B K_{ij}q_i^0q_j^0dv+\tau_T\, \int_B\frac{1}{K_m}\, \varphi_i\varphi_i dv\right]\right\}^{1/2},\label{4.20}
\end{split}
\end{equation}
then, from (\ref{4.19}) we obtain
\begin{equation}
\dot \Pi(t)-\sigma^2\Pi(t)\leq \frac{1}{\sqrt{2}}\, g(t),\label{4.21}
\end{equation}
which integrated provides the estimate (\ref{4.11}) and the proof is complete.

\medskip

\begin{remark}
We can establish estimates similar to (\ref{4.8}) and (\ref{4.11}) by starting with a conservation law in terms of the initial boundary value problem $\hat{\mathcal{P}}$.
\end{remark}

\begin{remark}
We have to outline that the estimate (\ref{4.8}) corresponds to the case of a stable system, while the estimate (\ref{4.11}) is of a new type in the sense that it allows the solutions to have an exponential growth in time and this should correspond to an unstable system (Cf. Quintanilla \cite{Quinta2002}).
\end{remark}

\section*{Spatial behavior results}

Throughout this section we study the influence of the delay times upon the spatial behavior of the transient and steady-state solutions within the context of the dual-phase-lag model of heat conduction.

\subsection*{Transient solutions}

In this subsection we assume that $B$ is the cylinder of length $L+h$ and of generic section $D$, more precisely we assume that $B=D\times (-h,L)$, $h>0$, $L>0$. Moreover, we suppose that the support of the given data $\mathcal{D}=\{r;T^0,q_i^0,\dot q_i^0;\vartheta,\xi\}$ is included into the closed cylinder $\overline{D}\times [-h,0]$. We are interested here into the behavior of the solution $\mathcal{S}=\{T,q_i\}$ of the initial boundary value problem $\mathcal{P}$ with respect to the distance $x_3$ from the support of loadings. To this aim we establish the following theorem of influence domain.

\begin{theorem}
Suppose that $k_{ij}$ is a positive definite tensor and assume that $a>0$. Moreover, we suppose that the delay times belong to the set defined by
\begin{equation}
\left\{0<\tau_q\leq 2\tau_T\right\}\cup \left\{0<2\tau_T<\tau_q\right\}.\label{5.1}
\end{equation}
Let $\mathcal{S}=\{T,q_i\}$ be a solution of the initial boundary value problem $\mathcal{P}$, in the cylinder $B=D\times (-h,L)$, $h>0$, $L>0$, corresponding to the given data $\mathcal{D}=\{r;T^0,q_i^0,\dot q_i^0;\vartheta,\xi\}$ having the support included into the closed cylinder $\overline{D}\times [-h,0]$. Then there exists a constant $c>0$, depending on the thermal constitutive coefficients and on the delay times, so that
\begin{equation}
\mathcal{S}(\mathbf{x},t)=\{T,q_i\}(\mathbf{x},t)=0\quad \textit{for all}\quad (\mathbf{x},t)\quad \textit{with}\quad x_3>ct,\quad t>0.\label{5.2}
\end{equation}

A possible value for $c$ is the following one
\begin{equation}
c_0=\sqrt{\frac{\kappa_M}{a_m}}\sqrt{\frac{2\tau_T}{\tau_q^2}+\frac{1}{\tau_T+\tau_q}}\label{5.3}
\end{equation}
when $0<\tau_q\leq 2\tau_T$; while when $0<2\tau_T<\tau_q$ a value for $c$ can be taken as
\begin{equation}
c_1=\frac{1}{\tau_q}\sqrt{\frac{\tau_T \kappa_M}{a_m}}\sqrt{\frac{5\tau_q^2+2\tau_T^2-6\tau_q\tau_T}{\tau_T^2+2\tau_q^2-3\tau_T\tau_q}},\label{5.4}
\end{equation}
where $a_m=\min_{\overline{B}}a$ and $\kappa_M=\max_{\overline{B}}k_M$.
\end{theorem}

\textbf{Proof.} First of all we note that
\begin{equation}
\hat{q}_i=-k_{ij}\left(T''_{,j}+\tau_T T'_{,j}\right),\quad \textit{for all}\quad \left(\mathbf{x},t\right)\in (D\times (0,L))\times (0,\infty),\label{5.5}
\end{equation}
so that, by means of the Cauchy-Schwarz inequality, we have
\begin{equation}
\begin{split}
&\hat{q}_i\hat{q}_i=-k_{ij}\hat{q}_i\left(T''_{,j}+\tau_T T'_{,j}\right)\leq \left(k_{ij}\hat{q}_i\hat{q}_j\right)^{\frac12}\left(k_{rs}T''_{,r}T''_{,s}\right)^{\frac12}\\[4mm]
&+\tau_T\left(k_{ij}\hat{q}_i\hat{q}_j\right)^{\frac12}\left(k_{rs}T'_{,r}T'_{,s}\right)^{\frac12}\leq \left(k_M\hat{q}_i\hat{q}_i\right)^{\frac12}\left[\left(k_{rs}T''_{,r}T''_{,s}\right)^{\frac12}+
\tau_T\left(k_{rs}T'_{,r}T'_{,s}\right)^{\frac12}\right].\label{5.6}
\end{split}
\end{equation}
Further, by using the arithmetic-geometric mean inequality, from (\ref{5.6}) we obtain
\begin{equation}
\hat{q}_i\hat{q}_i\leq k_M\left[\left(1+\varepsilon_1\right)k_{ij}T''_{,i}T''_{,j}+
\left(1+\frac{1}{\varepsilon_1}\right)\tau_T^2k_{ij}T'_{,i}T'_{,j}\right],\label{5.7}
\end{equation}
for all $\left(\mathbf{x},t\right)\in (D\times (0,L))\times (0,\infty)$ and for every positive parameter $\varepsilon_1$.

In order to study the spatial behavior of the solution $\mathcal{S}=\{T,q_i\}$ into the cylinder $\overline{D}\times [0,L]$ we introduce the following functional
\begin{equation}
I_{\delta}(x_3,t)=\int_0^t\int_0^s\int_{D_{x_3}}e^{-\delta z}\, \hat{T}(z)\hat{q}_3(z)dadzds,\quad x_3\in (0,L), \quad t\geq 0,\label{5.8}
\end{equation}
where $\delta \geq 0$ is a parameter whose values will be precisely given later. Further, we note that a direct differentiation with respect to the $x_3$ variable, and the use of the basic relations (\ref{2.6})-(\ref{2.9}), give
\begin{equation}
\begin{split}
&-\frac{\partial I_{\delta}}{\partial x_3}(x_3,t)=\frac12\, \int_0^t\int_{D_{x_3}}\, e^{-\delta s}\, \left[a \hat{T}^2(s)+\left(\tau_T+\tau_q+\delta \tau_q^2\right)k_{ij}T''_{,j}(s)T''_{,i}(s)\right.\\[4mm]
&\left.+\frac12\, \tau_T\tau_q^2 k_{ij}T'_{,j}(s)T'_{,i}(s)\right]dads+\frac{\tau_q^2}{4}\int_{D_{x_3}}\, e^{-\delta t}\, k_{ij}T''_{,j}(t)T''_{,i}(t)da\\[4mm]
&+\int_0^t\int_0^s\int_{D_{x_3}}\, e^{-\delta z}\, \left\{\frac{\delta a}{2}\, \hat{T}^2(z)+\left[1+\frac{\delta}{2}\left(\tau_T+\tau_q\right)+\frac{\delta^2}{4}\, \tau_q^2\right]k_{ij}T''_{,j}(z)T''_{,i}(z)\right.\\[4mm]
&\left.+\frac{1}{4}\, \tau_T\tau_q^2\left[\delta-2\left(\frac{1}{\tau_T}-\frac{2}{\tau_q}\right)\right]
k_{ij}T'_{,j}(z)T'_{,i}(z)\right\}dadzds,
\label{5.9}
\end{split}
\end{equation}
$x_3\in (0,L)$, $t\geq 0$.

At this time we can see that $-(\partial I_{\delta})/(\partial x_3) $ can be made positive if we choose for $\delta$ the value $\delta =\delta_0=0$ when $0<\tau_q\leq 2\tau_T$ and $\delta=\delta_1=2\sigma^2$ when $0<2\tau_T<\tau_q$. With these choices we deduce that
\begin{equation}
\begin{split}
&-\frac{\partial I_{\delta}}{\partial x_3}(x_3,t)\geq\frac12\, \int_0^t\int_{D_{x_3}}\, e^{-\delta s}\, \left[a \hat{T}^2(s)+\left(\tau_T+\tau_q+\delta \tau_q^2\right)k_{ij}T''_{,j}(s)T''_{,i}(s)\right.\\[4mm]
&\left.+\frac12\, \tau_T\tau_q^2 k_{ij}T'_{,j}(s)T'_{,i}(s)\right]dads\geq 0,\quad x_3\in (0,L), \, t\geq 0,
\label{5.10}
\end{split}
\end{equation}
where from now and in what follows we consider that $\delta$ has one of the two values $\delta_0=0$ and $\delta_1=2\sigma^2$.

Furthermore, we note that an integration in the relation (\ref{5.10}), with respect to $x_3$ variable on the interval $[x_3,L]$, gives
\begin{equation}
\begin{split}
&I_{\delta}(x_3,t)\geq\frac12\, \int_0^t\int_{\Omega_{x_3}}\, e^{-\delta s}\, \left[a \hat{T}^2(s)+\left(\tau_T+\tau_q+\delta \tau_q^2\right)k_{ij}T''_{,j}(s)T''_{,i}(s)\right.\\[4mm]
&\left.+\frac12\, \tau_T\tau_q^2 k_{ij}T'_{,j}(s)T'_{,i}(s)\right]dvds\geq 0,\quad x_3\in (0,L), \, t\geq 0,\label{5.11}
\end{split}
\end{equation}
where $\Omega_{x_3}=D\times [x_3,L]$. This proves that $I_{\delta}(x_3,t)$ can be considered as a measure of the solution $\mathcal{S}=\{T,q_i\}$.

On the other hand, by a direct differentiation with respect to time variable into relation (\ref{5.8}), we have
\begin{equation}
\frac{\partial I_{\delta}}{\partial t}(x_3,t)=\int_0^t\int_{D_{x_3}}\, e^{-\delta s}\, \hat{T}(s)\hat{q}_3(s)dads,\label{5.12}
\end{equation}
so that, by means of the Cauchy-Schwarz and the arithmetic-geometric mean inequalities and by using the estimate (\ref{5.7}), we obtain
\begin{equation}
\begin{split}
&\left| \frac{\partial I_{\delta}}{\partial t}(x_3,t) \right|\leq \frac12\, \int_0^t\int_{D_{x_3}}\, e^{-\delta s}\, \left[\varepsilon_2 a_m \hat{T}^2(s)+\frac{1}{\varepsilon_2 a_m}\, \hat{q}_3^2(s)\right]dads\\[4mm]
&\leq \frac12\, \int_0^t\int_{D_{x_3}}\, e^{-\delta s}\, \left\{\varepsilon_2 \left[a \hat{T}^2(s)\right]+\frac{\kappa_M\left(1+
\varepsilon_1\right)}{a_m\varepsilon_2\left(\tau_T+\tau_q+\delta\tau_q^2\right)}\right.\\[4mm]
&\left.\times \left[\left(\tau_T+\tau_q+\delta\tau_q^2\right)k_{ij}T''_{,i}(s)T''_{,j}(s)\right]+
\frac{2\tau_T\kappa_M\left(1+\varepsilon_1\right)}{a_m\varepsilon_1\varepsilon_2\tau_q^2}\left[\frac{\tau_T\tau_q^2}{2}\, k_{ij}T'_{,i}(s)T'_{,j}(s)\right]\right\}dads,\label{5.13}
\end{split}
\end{equation}
for all $\left(x_3,t\right)\in (0,L)\times (0,\infty)$ and for every positive parameters $\varepsilon_1$ and $\varepsilon_2$. At this stage we equate the coefficients of the various energy terms in (\ref{5.13}), that is we set
\begin{equation}
\varepsilon_2=\frac{\kappa_M\left(1+\varepsilon_1\right)}{a_m\varepsilon_2\left(\tau_T+
\tau_q+\delta\tau_q^2\right)}=\frac{2\tau_T\kappa_M\left(1+\varepsilon_1\right)}{a_m\varepsilon_1\varepsilon_2\tau_q^2}\label{5.14}
\end{equation}
and hence we fix
\begin{equation}
\varepsilon_1=\frac{2\tau_T\left(\tau_T+\tau_q+\delta\tau_q^2\right)}{\tau_q^2},\quad \varepsilon_2=\sqrt{\frac{\kappa_M}{a_m}}\sqrt{\frac{2\tau_T}{\tau_q^2}+\frac{1}{\tau_T+\tau_q+\delta\tau_q^2}}.
\label{5.15}
\end{equation}

With this choice, from the relations (\ref{5.10}) and (\ref{5.13}) we obtain the following first-order differential inequality
\begin{equation}
\left| \frac{\partial I_{\delta}}{\partial t}(x_3,t) \right|+c\, \frac{\partial I_{\delta}}{\partial x_3}(x_3,t)\leq 0,\quad \textit{for all}\quad \left(x_3,t\right)\in (0,L)\times (0,\infty),\label{5.16}
\end{equation}
where
\begin{equation}
c=\varepsilon_2=\sqrt{\frac{\kappa_M}{a_m}}\sqrt{\frac{2\tau_T}{\tau_q^2}+
\frac{1}{\tau_T+\tau_q+\delta\tau_q^2}}.\label{5.17}
\end{equation}

Further, the differential inequality (\ref{5.16}) can be treated as in Chiri\c t\u a and Quintanilla \cite{ChirQuinta1996} and Chiri\c t\u a and Ciarletta \cite{ChirCiar1999} in order to deduce that
\begin{equation}
I_{\delta}(x_3,t)=0\quad \textit{for all}\quad \left(x_3,t\right)\in (ct,L)\times (0,\infty),\label{5.18}
\end{equation}
and hence the relation (\ref{5.11}) implies the conclusion (\ref{5.2}) and the proof is complete.

\subsection*{Steady-state vibrations}

Throughout this subsection we assume that the cylinder $C=D\times (0,L)$ is made of a material with the delay times $\tau_q\geq 0$, $\tau_T\geq 0$ and that the conductivity tensor is positive definite. The cylinder is free of heat supply and it is thermally insulated on its lateral surface and on the end situated in the plane $x_3=L$. The cylinder is subjected to a harmonic perturbation on its base $x_3=0$ of the form
\begin{equation}
T(x_1,x_2,0,t)=h(x_1,x_2)e^{i\omega t},\quad (x_1,x_2)\in D_0,\quad t>0,\label{5.19}
\end{equation}
where $\omega >0$ is the frequency of perturbation and $i=\sqrt{-1}$ is the imaginary unit. Then inside of cylinder $C$ we will have the following harmonic vibration
\begin{equation}
\{T,q_r\}(\textbf{x},t)=\{\theta,Q_r\}(\textbf{x})e^{i\omega t},\label{5.20}
\end{equation}
where the amplitude $ \{\theta,Q_r\}$ of the vibration is a solution of the boundary value problem defined by the differential system
\begin{equation}
 \left(1+i\omega\tau_q-\frac12\, \tau_q^2\omega^2\right)Q_r=-\left(1+i\omega\tau_T\right)k_{rs}\theta_{,s},\label{5.21}
\end{equation}
\begin{equation}
 Q_{r,r}=-i\omega a\theta,\quad \textit{for all}\quad \mathbf{x}\in C,\label{5.22}
\end{equation}
with the boundary conditions
\begin{equation}
\theta(\mathbf{x})=0\quad \textit{on}\quad \left(\partial D_{x_3}\times (0,L)\right)\cup D_L,\label{5.23}
\end{equation}
and
\begin{equation}
\theta(x_1,x_2)=h(x_1,x_2),\quad (x_1,x_2)\in D_0.\label{5.24}
\end{equation}

We introduce the following functional
\begin{equation}
M(x_3)=-\int_{D_{x_3}}\left[\left(1+i\omega\tau_T\right)k_{3s}\theta_{,s}\overline{\theta}+\left(1-i\omega\tau_T\right)k_{3s}\overline{\theta}_{,s}\theta\right]da,\quad x_3>0,\label{5.25}
\end{equation}
where a superposed bar denotes complex conjugate. Further, we note that the relations (\ref{5.21}) and (\ref{5.25}) imply
\begin{equation}
M(x_3)=\int_{D_{x_3}}\left[\left(1+i\omega\tau_q-\frac12\, \tau_q^2\omega^2\right)Q_3\overline{\theta}+\left(1-i\omega\tau_q-\frac12\, \tau_q^2\omega^2\right)\overline{Q}_{3}\theta\right]da,\quad x_3>0,\label{5.26}
\end{equation}
and hence, by means of the divergence theorem and the use of relations (\ref{5.21}) and (\ref{5.23}), we get
\begin{equation}
-\frac{dM}{dx_3}(x_3)=2\int_{D_{x_3}}k_{rs}\theta_{,r}\overline{\theta}_{,s}da-2\tau_q\omega^2\, \int_{D_{x_3}}\, a\theta\overline{\theta} da.\label{5.27}
\end{equation}

On the other side, in view of the lateral boundary condition in (\ref{5.23}), we have
\begin{equation}
\int_{D_{x_3}}\theta_{,\alpha}\overline{\theta}_{,\alpha}da \geq \lambda \, \int_{D_{x_3}}\theta\overline{\theta}da,\label{5.28}
\end{equation}
where $\lambda$ is the lowest eigenvalue in the two-dimensional clamped membrane problem for the cross section $D_{x_3}$.

Now, if we use the estimate (\ref{5.28}) into relation (\ref{5.27}), we deduce that
\begin{equation}
-\frac{dM}{dx_3}(x_3)\geq 2\left(1-\frac{\tau_q a_M}{\lambda \kappa_m}\, \omega^2\right)\int_{D_{x_3}}k_{rs}\theta_{,r}\overline{\theta}_{,s}da,\label{5.29}
\end{equation}
where $a_M=\sup_{\overline{C}}|a| $ and $\kappa_m=\inf_{\overline{C}}k_m$. Further, we assume that the frequency of the vibration is such that
\begin{equation}
0< \omega < \sqrt{\frac{\lambda \kappa_m}{\tau_q a_M}},\label{5.30}
\end{equation}
and we note that an integration with respect to $x_3$ variable over $[x_3,L]$ and the use of the end boundary condition in (\ref{5.23}) give
\begin{equation}
M(x_3)\geq 2\left(1- \frac{\tau_q a_M}{\lambda \kappa_m}\omega^2\right)\int_{C_{x_3}}k_{rs}\theta_{,r}\overline{\theta}_{,s}dv\geq 0,\label{5.31}
\end{equation}
where $C_{x_3}=D\times (x_3,L)$.

By means of the Cauchy-Schwarz and the arithmetic-geometric mean inequalities and by use of the inequality (\ref{5.28}), from the relation (\ref{5.25}) we obtain
\begin{equation}
\left|M(x_3)\right|\leq \sqrt{\frac{\left(1+\tau_T^2\omega^2\right)}{\lambda \kappa_m^2}\, \sup_{\overline{C}}\left(k_{3s}k_{3s}\right)}\, \int_{D_{x_3}}k_{rs}\theta_{,r}\overline{\theta}_{,s}da.\label{5.32}
\end{equation}

Concluding, from the relations (\ref{5.30})-(\ref{5.32}), we obtain the first-order differential inequality
\begin{equation}
M(x_3)+\nu\, \frac{dM}{dx_3}\leq 0,\quad \textit{for all}\quad x_3\in (0,L),\label{5.33}
\end{equation}
where
\begin{equation}
\nu=\frac{\sqrt{\lambda}}{2\left(\lambda \kappa_m-\tau_q a_M\omega^2\right)}\, \sqrt{\left(1+\tau_T^2\omega^2\right) \sup_{\overline{C}}\left(k_{3s}k_{3s}\right)}.\label{5.34}
\end{equation}
When integrated, the differential inequality (\ref{5.33}) furnishes the estimate
\begin{equation}
0\leq M(x_3)\leq M(0)e^{-\frac{x_3}{\nu}}, \quad \textit{for all}\quad x_3\in (0,L),\label{5.35}
\end{equation}
that expresses the exponential decay of the amplitude $\{\theta,Q_r\}$ with respect to the distance $x_3$ at the loaded base.

Thus, we have the following theorem.
\begin{theorem}
Suppose that the hypothesis (H1) holds true. Then the spatial behavior of the amplitude $\{\theta,Q_r\}$ of the harmonic vibration (\ref{5.20}) is described by the inequality (\ref{5.35}), provided the frequency $\omega$ is lower than the critical value
 \begin{equation}
 \omega_c=\sqrt{\frac{\lambda \kappa_m}{\tau_q a_M}},\label{5.36}
 \end{equation}
 being the measure of the amplitude $M(x_3)$ as results from (\ref{5.31}), that is
  \begin{equation}
  M(x_3)\geq M^*(x_3),\quad M^*(x_3)=2\left(1-\frac{\omega^2}{\omega_c^2}\right)\int_{C_{x_3}}k_{rs}\theta_{,r}\overline{\theta}_{,s}dv.\label{5.37}
 \end{equation}
\end{theorem}

\section*{Results and Discussion}
In this paper we studied the dual-phase-lag model of heat conduction under various ranges of the delay times. The uniqueness and continuous data dependence problems are analyzed in the set $\{0\leq \tau_q\leq 2\tau_T\}\cup \{0<2\tau_T< \tau_q\}$. We established the two continuous dependence estimates (\ref{4.8}) and (\ref{4.11}) for any compact interval $[0,S]$, being the first one valid for $0\leq \tau_q\leq 2\tau_T $ and the second one for $0<2\tau_T< \tau_q$. We have to remark that the estimate (\ref{4.11}) allows the solutions to grow exponentially with respect to time variable. However, there is an open problem for the class of materials characterized by zero delay time of phase lag of the conductive temperature gradient and for which the delay time in the phase lag of heat flux vector is strictly positive. In such a case it should be expected to have an ill-posed model.

For the transient solutions we have established the theorem of influence domain as described by the relation (\ref{5.2}), where the speed of signal propagation is estimated by the value (\ref{5.3}) when the delay times satisfy the inequality $0< \tau_q\leq 2 \tau_T $ and by the value (\ref{5.4}) when the delay times satisfy the inequality $0<2\tau_T\leq \tau_q $. When $\tau_q=\tau_T=0$ the thermal model with delay times reduces to the classical model of heat conduction and this case can be treated following the method developed in Chiri\c t\u a and Ciarletta \cite{ChirCiar1999}, Chiri\c t\u a \cite{ChiritaJTS} and Quintanilla \cite{Quinta2001}. When $0=\tau_q<2\tau_T$ our above analysis fails to describe the spatial behavior of the transient solutions.

As regards to the steady-state solutions, we established an exponential decay estimate associated with the amplitude of vibration as described by (\ref{5.32}), provided the frequency of the vibration is lower than the critical value given by (\ref{5.33}) and the delay times satisfy $\tau_q\geq 0$ and $\tau_T\geq 0$. As it can be see the delay time $\tau_q$ influences the value of the critical frequency $\omega_c$, while $\tau_T$ and $\tau_q$ influence the decay rate of the amplitude vibrations.


\section*{References}



\begin{thebibliography}{99}

\bibliography{mybibfile}

\bibitem{Tzou1995a} D. Y. Tzou, A unified approach for heat conduction from macro to micro-scales, Journal of Heat Transfer 117 (1995) 8--16.

\bibitem{Tzou1995b} D. Y. Tzou, The generalized lagging response in small-scale and high-rate heating, International Journal of Heat and Mass Transfer 38 (1995) 3231--3234.

\bibitem{Tzou1995c} D. Y. Tzou, Experimental support for the lagging behavior in heat propagation, Journal of Thermophysics and Heat Transfer 9 (1995) 686--693.

\bibitem{Tzou2015} D. Y. Tzou, Macro- To Micro-Scale Heat Transfer:  The Lagging Behavior, John Wiley \& Sons, Chichester, 2015.

\bibitem{Chandra1998} D. S. Chandrasekharaiah, Hyperbolic thermoelasticity: A review of recent literature, Applied Mechanics Reviews 51 (1998) 705--729.

\bibitem{Fabrizio2014} M. Fabrizio, B. Lazzari, Stability and second law of thermodynamics in dual-phase-lag heat conduction, International Journal of Heat and Mass Transfer 74 (2014) 484--489.

\bibitem{Wang2001} L. Wang, M. Xu, X. Zhou, Well-posedness and solution structure of dual-phase-lagging
heat conduction, International Journal of Heat and Mass Transfer 44 (2001) 1659--1669.

\bibitem{Wang2002} L. Wang, M. Xu, Well-posedness of dual-phase-lagging heat equation: higher
dimensions, International Journal of Heat and Mass Transfer 45 (2002) 1165--1171.

\bibitem{Xu2002} M. Xu, L. Wang, Thermal oscillation and resonance in dual-phase-lagging heat
conduction, International Journal of Heat and Mass Transfer 45 (2002) 1055--1061.

\bibitem{Quinta2002} R. Quintanilla, Exponential stability in the dual-phase-lag heat conduction theory, Journal of Non-Equilibrium Thermodynamics 27 (2002) 217--227.

\bibitem{Horgan2005} C. O. Horgan, R. Quintanilla, Spatial behaviour of solutions of the dual-phase-lag heat
equation, Mathematical Methods in the Applied Sciences 28 (2005) 43--57.

\bibitem{Quinta2006} R. Quintanilla, R. Racke, A note on stability in dual-phase-lag heat conduction, International Journal of Heat and Mass Transfer 49 (2006) 1209--1213.

\bibitem{Racke2006} R. Quintanilla, R. Racke, Qualitative aspects in dual-phase-lag thermoelasticity, SIAM Journal on Applied Mathematics 66 (2006) 977--1001.

\bibitem{Quinta2007} R. Quintanilla, R. Racke, Qualitative aspects in dual-phase-lag heat conduction, Proceedings of the Royal Society of London A 463 (2007) 659--674.


\bibitem{Liu2007} K. C. Liu, P. C. Chang, Analysis of dual-phase-lag heat conduction in cylindrical system with a hybrid method, Applied Mathematical Modelling 31 (2007) 369--380.

\bibitem{Franchi2014} M. Fabrizio, F. Franchi, Delayed thermal models: Stability and thermodynamics, Journal of Thermal Stresses 37 (2014) 160--173.

\bibitem{Racke2015} R. Quintanilla, R. Racke, Spatial behavior in phase-lag heat conduction, Differential and Integral Equations 28 (2015) 291--308.

\bibitem{Abouelregal2015} A. E. Abouelregal, S. M. Abo-Dahab, Study of the dual-phase-lag model of thermoelasticity for a half-space problem with rigidly fixed surface in the presence of a thermal shock, Journal of Computational and Theoretical Nanoscience 12 (2015) 38--45.

\bibitem{Quinta2008} R. Quintanilla, A well posed problem for the dual-phase-lag heat conduction, Journal of Thermal
Stresses 31 (2008) 260--269.

\bibitem{Quinta2009} R. Quintanilla, A well-posed problem for the three-dual-phase-lag heat conduction, Journal of Thermal
Stresses 32 (2009) 1270--1278.

\bibitem{Quinta2016} R. Quintanilla, On uniqueness and stability for a thermoelastic theory, Mathematics and Mechanics of Solids DOI: 10.1177/1081286516634154.

\bibitem{Amendola2016} G. Amendola, M. Fabrizio, M. Golden, B. Lazzari, Second-order approximation for heat conduction: Dissipation principle and free energies, Proceedings of the Royal Society A: Mathematical, Physical and Engineering Sciences 472 (2016) art. no. 20150707.

\bibitem{Chirita2015} S. Chiri\c t\u a, M. Ciarletta, V. Tibullo, On the wave propagation in the time differential dual-phase-lag thermoelastic model, Proceedings of the Royal Society A: Mathematical, Physical and Engineering Sciences 471 (2015) art. no. 20150400.

\bibitem{ChirQuinta1996} S. Chiri\c t\u a, R. Quintanilla, On Saint-Venant's principle in linear elastodynamics. Journal of Elasticity 42 (1996) 201--215.

\bibitem{ChirCiar1999} S. Chiri\c t\u a, M. Ciarletta, Time-weighted surface power function method for the study of spatial behaviour in
dynamics of continua, European Journal of Mechanics - A/Solids 18 (1999) 915--933.

\bibitem{ChiritaJTS} S. Chiri\c t\u a, Spatial decay estimates for solutions describing harmonic vibrations in a thermoelastic cylinder, Journal of Thermal Stresses 18 (1995) 421--436.

\bibitem{Quinta2001} R. Quintanilla, End effects in thermoelasticity, Mathematical Methods in the Applied Sciences 24 (2001) 93--102.


\end{thebibliography}
\end{document}